\documentclass[reqno]{amsart}

\usepackage[utf8]{inputenc}
\usepackage[T1]{fontenc}

\usepackage{hyperref}
\usepackage{amsmath}
\usepackage{amssymb}
\usepackage{amsfonts}
\usepackage{graphicx}
\usepackage{amsthm}
\usepackage{enumerate}
\usepackage{lscape}
\usepackage{dsfont}
\usepackage{color}
\usepackage{mathtools}

\usepackage{setspace}
\newcommand{\R}{\mathds{R}}

\newcommand{\CP}{\mathds{C}\mathrm{P}}

\newcommand{\C}{\mathds{C}}            
\newcommand{\de}{\partial}          

\newcommand{\K}{K\"{a}hler}

\newcommand{\ov}[1]{\overline{#1}}
\newcommand{\wi}[1]{\widetilde{#1}}

\newcommand{\deb}{\ov\partial}

\newcommand{\ngh}{neighbourhood}
\newcommand{\lmb}{\lambda}

\newcommand{\W}{\Omega}
\newcommand{\w}{\omega}

\newcommand{\va}{\varphi}

\newcommand{\hs}{\hspace{0.1em}}

\newcommand{\Co}{\C^n\setminus \left\{0\right\}}
\newcommand{\tCh}{\tilde\C^n\setminus H}
\newcommand{\tC}{\tilde\C^n}

\def\F{\wi{\mathcal F}}

\newtheorem{theor}{Theorem}[section]
\newtheorem{prop}[theor]{Proposition}

\newtheorem{lem}[theor]{Lemma}

\begin{document}

\title{On holomorphic isometries  into blow-ups of $\C^n$}

\author{Andrea Loi}
\address{(Andrea Loi) Dipartimento di Matematica \\
         Universit\`a di Cagliari (Italy)}
         \email{loi@unica.it}

\author{Roberto Mossa}
\address{(Roberto Mossa) Dipartimento di Matematica \\
         Universit\`a di Cagliari (Italy)}
         \email{roberto.mossa@unica.it}

\thanks{
The first author  supported  by Prin 2015 -- Real and Complex Manifolds; Geometry, Topology and Harmonic Analysis -- Italy, by INdAM. GNSAGA - Gruppo Nazionale per le Strutture Algebriche, Geometriche e le loro Applicazioni.}


\begin{abstract}
We study  the  \K-Einstein manifolds which admits a  holomorphic isometry into either the generalized Burns-Simanca manifold $(\tilde \C^n, g_S)$ 
or the Eguchi-Hanson manifold $(\tilde \C^2, g_{EH})$.
Moreover, we prove that $(\tilde \C^n, g_S)$ and $(\tilde \C^2, g_{EH})$ are not relatives to any homogeneous bounded domain.
\end{abstract}
 
\maketitle

\tableofcontents

\section{Introduction}
The blow-up $\tilde \C^n$  of  $\C^n$  at the origin can be equipped with two  natural and interesting  \K\ metrics: the {\em Eguchi-Hanson metric} $g_{EH}$ on $\tilde \C^2$ and the
{\em generalized Burns-Simanca} metric $g_{S}$ on $\tilde \C^n$, $n\geq 2$. The \K\ form associated to the {\em Eguchi-Hanson}  metric
on $\C^2 \setminus \{0\}\cong \tilde{\C}^2 \setminus H$ given by
 \begin{equation}\label{eqpoteho}
\omega_{EH}=\frac{i}{2} \partial \bar{\partial}\left(\sqrt{\|z\|^4+1}+\log \|z\|^2-\log (1+\sqrt{\|z\|^4+1})\right), \ \   \|z\|^2=|z_1|^2+|z_2|^2,
\end{equation}
where $H\cong \C P^1$ denotes  the exceptional divisor arising by the blow-up construction (one can  show that  $\omega_{EH}$, a priori defined only on $\C^2 \setminus \{0\}$, extends to all $\tilde{\C}^2$).
The \K\   metric $g_{EH}$ is complete Ricci flat (not flat) (see \cite{EH}). The reader is also referred to \cite{Fraeh} and  \cite{LSZtwo} for other  properties  of this metric.
For $n\geq 2$ the generalized Burns-Simanca metric $g_{S}$ on $\tilde\C^n$ (first considered by  \cite{Simanca}) is the \K\ metric   whose \K\ form on $\C^n \setminus \{0\}\cong \tilde{\C}^n \setminus H$ is given by
\begin{equation}\label{eqpotsim0}
\omega_{S}=\frac{i}{2} \partial \bar{\partial}(\|z\|^2+\log \|z\|^2),\ \   \|z\|^2=|z_1|^2+\cdots+|z_n|^2,
\end{equation}
where  $H\cong \C P^{n-1}$ denotes  the exceptional divisor.
When $n=2$ the metric $g_{S}$ is simply called the {\em Burns-Simanca metric}
and it has been discovered by Burns \cite{Burns} when $n = 2$ and firstly described by Le Brun \cite{Brun}. 
This metric is scalar flat (but not Ricci-flat) and it is also an example (see \cite{Frasim}) of complete and non homogeneous metric admitting a regular quantization
(other properties of this metric related to the coefficients of TYZ expansion can be found in \cite{LSZthird}).
It is also worth mentioning that the Burns-Simanca metric has been  an important tool  in \cite{AP09} for the construction of cscK metrics on the blow-up of $\CP^2$
at a finite number of points.

In this paper we address two  problems.  The first one is about  the existence and uniqueness of 
 \K-Einstein (KE in the sequel) submanifolds of
either $(\tilde \C^n, g_S)$ or $(\tilde \C^2, g_{EH})$. 
The second one deals with the study of those \K\ manifolds which are relatives to either  $(\tilde \C^n, g_S)$ or $(\tilde \C^2, g_{EH})$.
Recall that in \cite{diloi}
the first author of the present paper together with A. J. Di Scala, inspired by Umehara's work \cite{Umeararel},  have christened  two \K\ manifolds $(S_1, g_1)$ and $(S_2, g_2)$ to be  {\em relatives} if they share a common (non trivial) \K\ submanifold, 
i.e. there exist a \K\ manifold $(M, g)$ of positive dimension  and two holomorphic isometries  $\varphi_j:M\rightarrow S_j$, $j=1, 2$ (the reader is referred to \cite{LMkimm} and references therein
for further results on relatives \K\ manifolds).

Our first result is the following theorem dealing  with the above mentioned  problems when the ambient space is the generalized Burns-Simanca metric.

\begin{theor}\label{mainteor}
Let $(M, g)$ be a complex $m$-dimensional ($m\geq 1$) \K\ manifold which admits a holomorphic isometry $\va :(M, g)\rightarrow (\tilde \C^n, g_S)$.
Then the following facts hold true.
\begin{itemize}
\item [(i)]
If $g$ is the flat metric then $m=1$,  $M$ is an open subset of  $\C$  and $\va$ is the restriction of 
the holomorphic isometry  $\Phi: (\C, g_0) \to (\tilde\C^n, g_S)$  given by
\begin{equation}
\label{eqpropflats}
\Phi(z)=\left((z  + \lmb )e , [e] \right),\ \lambda\in\C, \ e\in \C^n,\  \|e\|=1, 
\end{equation}
where $g_0$ denotes the flat metric.
\item [(ii)]
If  $g$ is KE with Einstein constant $\lambda \neq 0$ then 
 $\varphi (M)\subseteq H$, where $H$ is the exceptional divisor. 
Consequently,  if $n-m\leq 3$  then $M$ is either an open subset of the complex quadric or an open subset of a complex projective space totally geodesically embedded into $H\cong \C P^{n-1}$.
\item [(iii)]
Assume $n=2$. If $g$ is Ricci flat (i.e. KE with $\lambda =0$)  then $M$ is an open subset  of $\C$
and $\va$ is given by the restriction of 
\eqref{eqpropflats} to $M$.
\item [(iv)]
any \K\ submanifold of $(M, g)$ does not admit a holomorphic isometry  into a homogeneous bounded domain.
\end{itemize}
\end{theor}

Notice that (i) is a sort of rigidity result for flat submanifolds of $(\tilde \C^n, g_S)$.
Regarding (ii) observe that,  despite the fact that a  KE submanifold  of a non elliptic complex space form  is forced to be  totally geodesic \cite{UmearaE},  the classification of the KE submanifolds of the complex projective space is still missing (the reader is referred to \cite[Chapter 3]{LoiZedda-book} for an updated material on this subject).
In \cite{LZ} (see also \cite{LSZ} and \cite{LSZpac}) the authors have conjectured  that these manifolds are indeed open subset of flag manifolds (i.e.  compact and simply-connected homogeneous \K\ manifolds).
Therefore, in view of (ii) we believe that a KE manifold admitting a holomorphic isometry into   $(\tilde \C^n, g_S)$ is holomorphically isometric to an open subset of a flag manifold.
About point (iii) we still do not know if there exist Ricci flat (not flat) \K\ submanifolds of  $(\tilde \C^n, g_S)$ for $n\geq 3$.
Finally (iv) shows that the generalized Burns-Simanca manifold is not relative to any bounded homogeneous domain.  

Our second and last result   shows the analogous of (ii) and (iv)  for the Eguchi-Hanson metric.

\begin{theor}\label{mainteor2}
Let $(M, g)$ be a complex $m$-dimensional ($m\geq 1$) \K\ manifold which admits a holomorphic isometry into  $(\tilde \C^2, g_{EH})$.
Then the following facts hold true.
\begin{itemize}
\item [(a)]
If  $g$ is KE with Einstein constant $\lambda \neq 0$ then  $m=1$ and 
 $\varphi (M)\subseteq H$.
\item [(b)]
$(M, g)$ does not admit a holomorphic isometry  into a homogeneous bounded domain.
\end{itemize}
\end{theor}

We still do not know if $(\tilde \C^2, g_{EH})$ admits a  flat \K\   submanifold 
similarly to that of $(\tilde\C^n, g_S)$ given by $\Phi (\C)$ where  $\Phi$ is the map \eqref{eqpropflats}. Observe that  (v) in Theorem \ref{mainteor} and (b) in Theorem \ref{mainteor2} can be considered  extensions of  \cite[(ii) Theorem 1.1]{LMkimm} where the authors of the present paper shows that the definite or indefinite flat space is not relative to any  bounded homogenous domain.

 The proofs  of Theorem \ref{mainteor} and Theorem \ref{mainteor2} can be found in Section \ref{proofs}.
The main tools in their proofs are a detailed analysis of Calabi's diastasis function of the generalized Burns-Simanca metric and Eguchi-Hanson metric
(Proposition \ref{propfond} in  the next section)
and some trascendental  properties of holomorphic Nash algebraic function 
 (Lemma \eqref{lemnash}) already used in \cite{LMpams} and \cite{LMkimm}.

\section{Calabi's diastasis function of $g_S$ and $g_{EH}$}\label{necessary}
The blow-up $\tilde \C^n$ of $\C^n$ at the origin is the closed submanifold of $ \mathbb{C}^{n} \times \mathbb{C} P^{n-1}$ given by
$$
\tilde{\mathbb{C}}^{n}:=\left\{(z,[t]) \in \mathbb{C}^{n} \times \mathbb{C} P^{n-1} \mid t_{\alpha} z_{\beta}-t_{\beta} z_{\alpha}=0,\ 1 \leq \alpha < \beta \leq n\right\}.
$$ 
Let $U_{j}:=\left\{\left[z_{1}: \ldots: z_{n}\right] \mid z_{j} \neq 0\right\} \subset \mathbb{C} P^{n-1}$ and $\tilde{U}_{j}=\left(\mathbb{C}^{n} \times U_{j}\right) \cap \tilde{\mathbb{C}}^{n}$. A complex atlas  for $\tilde{\mathbb{C}}^{n}$ is given by $(\tilde U_j, \psi_j)$, where  $\psi_{j}: \tilde{U}_{j} \rightarrow \mathbb{C}^{n}$,
$$
\psi_{j}\left(\left(z_{1}, \ldots, z_{n}\right),\, \left[t_{1}: \ldots: t_{n}\right]\right) =\left(\frac{t_{1}}{t_{j}}, \ldots, \frac{t_{j-1}}{t_{j}}, z_{j}, \frac{t_{j+1}}{t_{j}}, \ldots, \frac{t_{n}}{t_{j}}\right), \ j=1,\dots,n, $$ 
whose inverse is
 \begin{equation}\label{eqloccoor}
\psi_{j}^{-1}(w)= \left(\left(w_{j} w_{1}, \ldots, w_{j} w_{j-1}, w_{j}, w_{j} w_{j+1}, \ldots, w_{j} w_{n}\right),\ \left[w_{1}: \ldots: w_{j-1}: 1: w_{j+1}: \ldots: w_{n}\right]\right),
 \end{equation}
 with $w=\left(w_{1}, \ldots, w_{n}\right)$. 
 If  $H=\left\{(0,[t]) \in \tilde{\mathbb{C}}^{n}\right\} \simeq \mathbb{C} P^{n-1}$ denotes  the exceptional divisor, then the map 
\begin{equation}\label{eqcoordsetm}
p_r:\tilde{\mathbb{C}}^{n} \setminus H \to \C^n\setminus \left\{0\right\}, \qquad (z,[t]) \mapsto z
\end{equation}
is a biholomorphism with inverse $p_r^{-1}(z)=\left(z,[z]\right)$. 

In the proof of our main results we need the following proposition where we describe Calabi's diastasis functions of   the  generalized Burns-Simanca metric   and the Eguchi-Hanson metric and we verify   that the restrictions of these metrics to the exceptional  divisor has 
constant positive holomorphic sectional curvature.
Recall that given a    real analytic \K\  metric $g$ on a complex manifold manifold $M$ and a point $p\in M$, 
{\em Calabi's diastasis function}  $D^g_p:U\to\R$  at $p$ is defined as
\begin{equation}\label{defdiastasis}
D^g_p(z)=\tilde\psi(z,\bar z)+ \tilde\psi(p,\bar p) -\tilde\psi(z,\bar p)-\tilde\psi(p,\bar z),
\end{equation}
where $\tilde{\psi}: U  \times U\rightarrow \R$ is a complex analytic extension 
obtained by duplicating the variables $z$ and $\bar{z}$ of 
a \K\ potential  $\psi:U\rightarrow \R$ for the metric  $g$ and where $U$ is a neighborhood of $p$.
\begin{prop}\label{propfond}
Let $\tilde \C^n$  be the blow-up of $\C^n$ at the origin. Then:
\begin{itemize}
\item [1.]
Calabi's diastasis function of $g_S$ centered at $(q,[t])\in \tCh$, with respect  to the coordinates \eqref{eqcoordsetm},  is given  by
\begin{equation}\label{eqpotsim}
D^{g_S}_{(q,[t])}(z)=\left\|z-q\right\|^2 +  \log \frac{\left\|z\right\|^2\left\|q\right\|^2}{\left|z\cdot \ov{q}\right|^2}.
\end{equation}
\item [2.]
Calabi's diastasis function  of $g_{EH}$
centered at $(q,[t])\in \tilde \C^2$
is given by
{\small
\begin{equation}\begin{split}\label{eqpoteh} 
D^{g_{EH}}_{(q, [t])}(z)&=\sqrt{\|z\|^{4}+1}+ \sqrt{\|q\|^{4}+1}-
\sqrt{\left(z\cdot\ov{q}\right)^{2}+1}-
\sqrt{\left(\ov{z}\cdot{q}\right)^{2}+1}\\
 & + \frac{\|z\|^{2}\|q\|^{2}\left|1+\sqrt{\left(z\cdot\ov{q}\right)^2+1}\right|^2}{|z\cdot \ov{q}|^2\left(1+\sqrt{\|z\|^{4}+1}\right)\left(1+\sqrt{\|q\|^{4}+1}\right)}
\end{split}\end{equation}}
\end{itemize}
Moreover, 
$\left(H,{g_{{S}}}_{\big | H}\right)$ (resp. $\left(H,{g_{{EH}}}_{\big | H}\right)$) is holomorphically isometric to the complex projective space $\left(\C P^{n-1}, g_{FS}\right)$ (resp. $\left(\C P^{1}, g_{FS}\right)$) equipped with the Fubini-Study metric $g_{FS}$ of holomorphic sectional curvature $4$.
\end{prop}
\proof
Expressions \eqref{eqpotsim} and \eqref{eqpoteh} follow by combining  \eqref{eqpotsim0} and \eqref{eqpoteho} with the definition of Calabi's diastasis function.

In local coordinates $(\tilde U_j,\psi_j)$, the metric $g_S$ has a  \K\ potential given by
\begin{equation}\label{eqlocpottuj}
\phi^{S}_j(w)=\left(\left|w_{j}\right|^{2}\left(1+\|w\|^{2}-\left|w_{j}\right|^{2}\right)+\log \left(1+\|w\|^{2}-\left|w_{j}\right|^{2}\right)\right), \qquad j=1,\dots,n.
\end{equation}
Thus the restriction of $g_S$ to the exceptional divisor reads as 
\begin{equation}\label{eqsimancafubini}
{\phi^{ S}_j}_{\big |  \tilde U_j\cap H}(w)=\phi^{S}_j(w_1,\dots,w_{j-1},0,w_{j+1},\dots,w_n)=\log \left(1+\left\|w\right\|^{2}\right),
\end{equation}
which is a \K\ potential of the Fubini-Study metric of $H=\CP^{n-1}$ of holomorphic sectional curvature $4$.

Similarly, in local coordinates $(\tilde U_j,\psi_j)$, $j=1,2$ given by  \eqref{eqloccoor}),  the associated \K\ potential reads 
$$
\phi^{EH}_{\tilde U_1}(w)=\sqrt{\left|w_{1}\right|^{4}\left(1+\left|w_{2}\right|^{2}\right)^{2}+1}+\log \left(\frac{1+\left| w_{2}\right|^{2}}{1+\sqrt{\left|w_{1}\right|^{4}\left(1+\left|w_{2}\right|^{2}\right)^{2}+1}}\right)
$$
$$
\phi^{{EH}}_{\tilde U_2}(w)=\sqrt{\left|w_{2}\right|^{4}\left(1+\left|w_{1}\right|^{2}\right)^{2}+1}+\log \left(\frac{1+\left| w_{1}\right|^{2}}{1+\sqrt{\left|w_{2}\right|^{4}\left(1+\left|w_{1}\right|^{2}\right)^{2}+1}}\right),
$$
respectively.
By restricting them to the exceptional divisor, we get
\begin{equation}\label{eqrestrhgeh}
{\phi^{EH}_{\tilde U_j}}_{\big |  \tilde U_j\cap H}(w)=1+\log \left(\frac{1+\left| w_{j}\right|^{2}}{2}\right) \qquad j=1,2,
\end{equation}
and so 
$\left(H,{g_{{EH}}}_{\big | H}\right)$ is holomorphically isometric to the complex projective space $\left(\C P^1, g_{FS}\right)$ equipped with the Fubini-Study metric $g_{FS}$ of holomorphic sectional curvature $4$.
\endproof

\section{Proof of the main results}\label{proofs}

Let $\mathcal{N}^{m}$ be the set of real analytic functions $\xi: V \subset \mathbb{C}^{m} \rightarrow \mathbb{R}$ defined in some open neighbourhood $V \subset \mathbb{C}^{m}$, such that its real analytic extension $\tilde{\xi}(z, w)$ in a neighbourhood of the diagonal of $V \times \operatorname{Conj} V$ is a holomorphic Nash algebraic function (for background material on Nash functions, we refer the readers to \cite{XJHuang} and \cite{Tworz}). We define
\begin{equation*}
\mathcal{F}=\left\{\xi\left(f_{1}, \ldots, f_{m}\right) \mid \xi \in \mathcal{N}^{m}, f_{j} \in \mathcal{O}_{0}, j=1, \ldots, m, m>0\right\}
\end{equation*}
where $\mathcal{O}_{0}$ denotes the germ of holomorphic functions around $0 \in \mathbb{C}$ and we set
\begin{equation*}
\tilde{\mathcal{F}}=\{\psi \in \mathcal{F} \mid \psi \text { is of diastasis-type }\}
\end{equation*}
Here we say (see also \cite{LMpams}) that a real analytic function defined on a neighborhood $U$ of a point $p$ of a complex manifold $M$ is of diastasis-type if in one (and hence any) coordinate system $\left\{z_{1}, \ldots, z_{n}\right\}$ centered at $p$ its expansion in $z$ and $\bar{z}$ does not contain non constant purely holomorphic or anti-holomorphic terms (i.e. of the form $z^{j}$ or $\bar{z}^{j}$ with $j>0$ ). Clearly the diastasis $D^g_p$ is  a function of diastasis-type.

In the proofs of our main results we need the following lemmata.

\begin{lem}\label{lemnash}(\cite[Theorem 2.1]{LMkimm})
Let $\psi_{0} \in \widetilde{\mathcal{F}} \backslash \mathbb{R}$. Then for every $\mu_{1}, \ldots, \mu_{\ell} \in \mathbb{R}$ we have
\begin{equation*}
e^{\psi_{0}} \notin \tilde{\mathcal{F}}^{\mu_{1}} \cdots \tilde{\mathcal{F}}^{\mu_{\ell}} \backslash \mathbb{R}
\end{equation*}
where $\widetilde{\mathcal{F}}^{\mu_{1}} \ldots \widetilde{\mathcal{F}}^{\mu_{\ell}}=\left\{\psi_{1}^{\mu_{1}} \cdots \psi_{\ell}^{\mu_{\ell}} \mid \psi_{1}, \ldots, \psi_{\ell} \in \widetilde{\mathcal{F}}\right\}$
\end{lem}

\begin{lem}\label{thmsimkegen}
Let $(N, h)$ be a \K\ manifold and let $\left\{z_1,\dots,z_n\right\}$ be a system of coordinates for $N$ around $q\in N$. Assume that the diastasis $D_q^{h}$ associated to $h$ and centred in $q$ is of the form
$$
D_q^{h} (z) = \Psi(z) + \log \left(\psi_1^{\gamma_1}(z)\cdots \psi_r^{\gamma_r}(z)\right)
$$
were $\psi_1,\dots,\psi_r\in \F$, $\gamma_1,\dots,\gamma_r\in\R$ and $\Psi\in\F$ is strictly plurisubharmonic at $q$.
Then any  KE submanifold $\left(M,g \right)$ of $(N, h)$ passing through $q$, is Ricci flat.
\end{lem}
\proof
Assume by contradiction that $f:(M,g) \to (N, h)$ is a holomorphic isometry passing through $q$  and that $g$ is KE non Ricci flat. Let $q=f(p)$. Fixed a system of coordinates $\left\{w_1,\dots,w_m\right\}$ for $M$  centered at $p$, the diastasis associated to $g_M$ and centered at $p$ is given by:
\begin{equation}\label{eqdgn=dq}
\end{equation}
$$
D^{g}_p(w)=D_q^{h} (f(w))=\Psi(f(w)) + \log \left(\psi_1^{\gamma_1}(f(w))\cdots \psi_r^{\gamma_r}(f(w))\right).
$$
From the uniqueness of the diastasis function, we see that the Einstein condition
$
-i\de\deb \log \det \left(\de\deb D^{g}_p \right)= \lmb\frac{i}{2}\de\deb D^{g}_p
$
is equivalent to the following equation (see e.g. \cite[Proof of Proposition 4.1]{LMkimm})
$$
\det \left(\de\deb D^{g}_p \right)= e ^{-\frac \lmb 2D^{g}_p},
$$
By substituting \eqref{eqdgn=dq} in the previous equation we get
$$
\det \left(\de\deb D^{g}_p \right)=e^{-\frac \lmb 2\Psi(f(w))}\hs\psi_1^{-\frac {\lmb \gamma_1} {2}}(f(w))\cdots \psi_r^{-\frac {\lmb \gamma_r} {2}}(f(w)).
$$
It is not hard to see that $\det \left(\de\deb D^{g}_p \right)\in\F$ (see \cite[Proof of Proposition 4.1]{LMkimm}), so that
$$
e^{{-\frac \lmb 2\Psi(f(w))}} 
\in \F^{\frac {\lmb \gamma_1} {2}}\cdots \F^{\frac {\lmb \gamma_r} {2}}\F.
$$
By Lemma \ref{lemnash}, we see that $\Psi(f(w))$ is forced to be zero,
By hypothesis, $\Psi$ is strictly plurisubharmonic at $q$, in particular  $\Psi$ must be a diastasis function for a \K\ metric $g_\Psi$ defined in a \ngh\ of $q$. By \cite[Proposition 5]{Cal} we know that if $\rho({q}, {q_0})$ is the geodesic distance induced by $g_\Psi$, between $q$ and $q_0$, then 
$$
\Psi (q_0)=\rho({q}, {q_0})^{2}+O\!\left(\rho({q}, {q_0})^{4}\right),
$$ 
in a \ngh\ of $q$. Therefore,  \eqref{eqpsifw=} implies that $f(w)=f(p)$ 
on a \ngh\ of $p$. In particular $f$ can not be a \K\ immersion,  contradicting the hypothesis. The proof is complete. 
\endproof
\begin{proof}[Proof of Theorem \eqref{mainteor}]
{\em Proof of (i)}
Notice that  the generalized Simanca metric $g_S$ is projectively induced, i.e. there exists a holomorphic isometry $\psi :M\rightarrow \C P^{\infty}$ into the infinite dimensional complex projective space equipped with the Fubini-Study metric (see \cite{Frasim} for a proof).  Thus by \cite[Theorem 13]{Cal}  $M$
is forced to be  an open
subset $V\subset \C^m$ equipped with the flat metric.
Hence one needs to show that if  $\va:(V,g_0)\to (\tilde \C^n,g_{ S })$ is a holomorphic isometry then $m=1$ and $\va$ is the restriction of the map \eqref{eqpropflats}.
We claim that $\va(V)$ pass through a point $q \in \tCh$. Assume by contradiction that $\va(V)\subset H$. From \eqref{eqsimancafubini} we know that the metric induced by $g_S$ on $H\equiv \C P^{n-1}$ is the Fubini-Study metric metric $g_{FS}$ hence $g_0=\va^*g_{FS}$
in constrast with Calabi's result \cite{Cal}, proving our claim.  

It is not restrictive to assume that $V$ is a \ngh\ of the origin of $\C^m$ and that  $\va(0)=q$.
From Calabi diastasis's hereditary property, in a \ngh\ $W\subset V$ of the origin, we have $D_0^{g_0}(z)=D^{g_{ S }}_{\va(0)}(\va(z))$. Assume also $ \va\left(W\right) \subset \tCh$. From \eqref{eqpotsim},  we deduce that 
\begin{equation*}
\|z\|^2= \left\|\va(z)-\va(0)\right\|^2 +  \log \frac{\left\|\va(z)\right\|^2\left\|\va(0)\right\|^2}{\left|\va(z)\cdot \ov{\va(0)}\right|^2}, \qquad z\in W.
\end{equation*} 

We can apply Lemma \ref{lemnash} to conclude that 
\begin{equation}\label{eqpropflatnash}
\|z\|^2- \left\|\va(z)-\va(0)\right\|^2=0 \quad\text{ and }\quad \frac{\left\|\va(z)\right\|^2\left\|\va(0)\right\|^2}{\left|\va(z)\cdot \ov{\va(0)}\right|^2}=1.
\end{equation}
Since the isometric action of $U(n)$ on $\tC$ (given by $U\cdot (z,[t])= \left(U\cdot z, [U\cdot t]\right)$) is the standard isometric action  on $\left(\C^n\setminus \left\{0\right\}, g_S\right)$  we can assume that $\va(0)=(0,\dots,0,\lmb)$, $\lmb\in\C^*$. By  the second equality in \eqref{eqpropflatnash}, we see that $\va(z)$ is proportional to $\va(0)$, that is
\begin{equation}\label{eqproorelsimfl}
\va(z)=(0,\dots,0,\va_n(z)), \qquad z\in W,
\end{equation}
where $\va_n: W\rightarrow\C$ is a holomorphic function.
Since $\va$ is an immersion, the previous equation proves  that $\dim W=m=1$. 
Moreover,  we see that in  the coordinates charts  \eqref{eqloccoor},  $\va(W)\subset\tilde U_n$ and $\va(W)\cap\tilde U_j=\emptyset$, $j=1,\dots,n-1$.

Since,  with respect to the coordinates $\psi_n$ one has
$$
\va(z) = \psi_n \circ p_r^{-1}\left(0,\dots,0,\va_n(z)\right)= \left(0,\dots,0,\va_n(z)\right), 
$$
equality $D_0^{g_0}(z)=D^{g_{ S }}_{\va(0)}(\va(z))$ yields
$$|z|^2=|\va_n(z) -\lmb|^2, \ z\in W \subset \C.$$
Deriving this relation with respect $\ov z$, we get
$$
z=\va_n (z)\ov{\va_n^\prime (z)}-{\lmb}\ov{\va_n^\prime (z)}=\ov{\va_n^\prime (z)}\left({\va_n (z)}-{\lmb}\right)
$$
yielding $\ov{\va_n^\prime (z)}=\mu\in \C$  and  $\va_n (z)={\lmb}+\mu^{-1} z$. 
On the other hand if  $\w_S$ and $\w_0$ are the \K\ forms associated to  $g_S$ and $g_0$, from \eqref{eqlocpottuj} and \eqref{eqproorelsimfl} we see that
$$
\frac{i}{2}\de\deb |z|^2=\w_0=\va^*\w_S=\frac{i}{2}\de\deb\left|{\lmb}+\mu^{-1} z\right|^{2}.
$$
Hence $\mu=e^{i\theta}$ which, up to the  $U(n)$-action on $\tC$ we can assume $\mu=1$. From the analyticity of $\va$, we see that $\va(z)=(0, \dots, 0, \lambda+z)$  for every $z\in V$, i.e. 
$\va =\Phi_{|V}$. 
{\em Proof of (ii)}
Fix $q\in \tCh\equiv\Co$. 
Since  $|z- q|^2$ is  strictly plurisubharmonic at $q$  the diastasis of $g_S$ given by \eqref{eqpotsim} satisfies the hypothesis of Lemma \ref{thmsimkegen}.  
We conclude that if $\lmb\neq 0$, then $\va (M) \subset H \equiv \C P^{n-1}$. If $\va (M)$ is open in $H$,  by \eqref{eqsimancafubini} we see that $g$ is the Fubini-Study metric. While, if $\va (M)$ is a submanifold of $H$ of codimension $1$ or $2$, due to S.-S. Chern \cite{chern} and K. Tsukada \cite{tsukada} results respectively, $M$ is either an open subset of the complex quadric or an open subset of a complex projective space holomorphically and totally geodesically embedded in $H$. 

{\em Proof of (iii)}
It follows by combining  (i) to the fact that a  Ricci flat metric on a complex one-dimensional manifold 
is flat.

{\em Proof of (iv)}
We need to show that it can not exist a \ngh\ $V\subset \C$ of the origin and maps  $\eta: V \to M$ and $\zeta :V\rightarrow \Omega$ such that
\begin{equation}\label{eqfondrel}
\eta^*g= \zeta^*g_\W .
\end{equation}
 Let us suppose by contradiction that such maps exist. We claim that  $\va(\eta(V))$ is not contained in $H$. Indeed from \eqref{eqsimancafubini} we know that the metric induced by $g_S$ on $H\equiv \C P^{n-1}$ is the Fubini-Study metric $g_{FS}$, hence $\eta^*(\va^*g_{FS})=\zeta^*g_{\Omega}$ in contrast with  \cite[Theorem 1.2]{Mossa}, claming that a projective manifold is not relative to a homogeneous bounded domain. Hence, we can  assume that $\va(\eta(0))=q \in \tCh$.
 
From \cite[Theorem 3.1]{LMkimm} (see also \cite{LMdedicata}), we now that there exists rational (and hence Nash) functions $F_{1}, \ldots, F_{s}$ and positive numbers $\delta_{1}, \ldots, \delta_{s}$ such that $\sum_{k=1}^{s} \delta_{k} \log F_{k}(z)$ is a \K\ potential of $g_\W$. Hence, by \eqref{eqpotsim} in a \ngh\ of $0$, \eqref{eqfondrel} yields  the following equation in terms of diastasis functions
\begin{equation}\label{eqinter}
\|\va(\eta (z))-q\|^2 = \log\left[ \prod_{k=1}^{s} \left(\frac{F_{k}\left(\zeta(z), \overline{\zeta(z)}\right) F_{k}\left(\zeta(0), \overline{\zeta(0)}\right)}{F_{k}\left(\zeta(z), \overline{\zeta(0)}\right) F_{k}\left(\zeta(0), \overline{\zeta(z)}\right)}\right)^{\delta_k}\left(\frac{\left\|\va(\eta (z))\right\|^2\left\|q\right\|^2}{\left|\va(\eta (z))\cdot \ov{q}\right|^2}\right)^{-1}\right]
\end{equation}
By Lemma \ref{lemnash}, we see that $\|\va(\eta (z))-q\|^2$  is forced to be a constant. Since  $\va(\eta (0))=q$ we see that $\va(\eta (z))=q$
for all  $z$ in a \ngh\ of $0$.   In particular $\eta$ cannot be an immersion, contradicting the hypothesis. 
\end{proof}

\begin{proof}[Proof of Theorem \eqref{mainteor2}]
{\em Proof of (a)}
Let  $q\in \tilde\C^2\setminus H$.
By a straightforward computation, we can see that 
\begin{equation}\label{eqnew}
\sqrt{|z|^{4}+1}+ \sqrt{|q|^{4}+1}-
\sqrt{\left(z\cdot\ov{q}\right)^{2}+1}
-\sqrt{\left(\ov{z}\cdot{q}\right)^{2}+1}
-\log \left(1+\sqrt{|z|^{4}+1}\right)
\end{equation}
is strictly plurisubharmonic at $q$. Hence, the diastasis of   $g_{EH}$ given by \eqref{eqpoteh} satisfies the hypothesis of Lemma \ref{thmsimkegen}.  
We conclude that if $\lmb\neq 0$, then $\va (M) \subset H \equiv \C P^{1}$.
{\em Proof of (b)}
The proof follows the same line of the proof of (iv) in Theorem \ref{mainteor}
by using  in \eqref{eqinter} the function \eqref{eqnew} instead of $\|\va (z)-q\|^2$.
\end{proof}

\end{document}